\documentclass[12pt]{amsart}

\usepackage{amsmath,amsthm,amssymb}

\font\teneufm=eufm10 \font\seveneufm=eufm7 \font\fiveeufm=eufm5
\newfam\frakturfam
\textfont\frakturfam=\teneufm \scriptfont\frakturfam=\seveneufm
\scriptscriptfont\frakturfam=\fiveeufm

\newtheorem{lm}{Lemma}
\newtheorem{theor}{Theorem}
\newtheorem{co}{Corollary}
\newtheorem{rem}{Remark}
\newtheorem{prob}{Problem}
\def\bee{\begin{eqnarray}}
\def\bes{\begin{eqnarray*}}
\def\eee{\end{eqnarray}}
\def\ees{\end{eqnarray*}}

\def\a{\alpha}

\def\Proof{{\sl Proof.}\ }

\title{The freiheitssatz and the automorphisms of free right-symmetric algebras}

\begin{document}
\date{}
\maketitle

\begin{center}

{\bf Daniyar Kozybaev}\\
Department of Mathematics, Eurasian National University\\
 Astana, 010008, Kazakhstan \\
e-mail: {\em kozybayev@rambler.ru},\\

{\textbf{Leonid Makar-Limanov} \footnote{Supported
by an NSA grant.}}\\
Department of Mathematics \& Computer Science,\\
Bar-Ilan University, 52900 Ramat-Gan, Israel and\\
Department of Mathematics, Wayne State University, \\
Detroit, MI 48202, USA\\
e-mail: {\em lml@math.wayne.edu}, \\

and\\

{\bf Ualbai Umirbaev}\\
Department of Mathematics, Eurasian National University\\
 Astana, 010008, Kazakhstan \\
e-mail: {\em umirbaev@yahoo.com}\\

\end{center}

\begin{abstract}
We prove the Freiheitssatz for right-symmetric algebras and the
decidability of the word problem for right-symmetric algebras with
a single defining relation. We also prove that two generated
subalgebras of free right-symmetric algebras are free and
automorphisms of two generated free right-symmetric algebras are
tame.
\end{abstract}


\noindent


\section{Introduction}

\hspace*{\parindent}
\bigskip

The Freiheitssatz (``freedom/independence theorem" in German), one
of the most important theorems of combinatorial group theory, was
proposed by M.\,Dehn in the geometric setting and proved by his
student, W.\, Magnus, in his doctoral thesis \cite{Magnus}. The
Freiheitssatz says the following: Let $G=\langle x_1,x_2,\ldots,
x_n | r=1\rangle$ be a group defined by a single cyclically
reduced relator $r$. If $x_n$ appears in $r$, then the subgroup of
$G$ generated by $x_1,\ldots, x_{n-1}$ is a free group, freely
generated by $x_1,\ldots, x_{n-1}$. W.\,Magnus also proved in
\cite{Magnus} the decidability of the word problem for groups with
a single defining relation. The Freiheitssatz for solvable and
nilpotent groups was studied by N.\,S.\,Romanovskii
\cite{Romanovskii}.

The Freiheitssatz and the decidability of the word problem for Lie
algebras with a single defining relation was proved by
A.\,I.\,Shirshov \cite{Shir2}. L.\,Makar-Limanov \cite{Makar2}
proved the Freiheitssatz for associative algebras over a field of
characteristic zero. Question about the decidability of the word
problem for associative algebras (and also for semigroups) with a
single defining relation and the Freiheitssatz for associative
algebras in positive characteristic remain open (see
\cite{BokutKukin}).

An algebra $A$ over an arbitrary field $k$ is called {\em
right-symmetric} if it satisfies the identity \bee\label{f1}
(xy)z-x(yz)=(xz)y-x(zy). \eee In other words, the associator
$(x,y,z)=(xy)z-x(yz)$ is symmetric in $y$ and $z$. The variety of
right-symmetric algebras is Lie-admissible, i.e., each
right-symmetric algebra $A$ with the operation $[x,y]=xy-yx$ is a
Lie algebra. Right-symmetric algebras are associated with locally
affine manifolds (see \cite{Vinberg}).

A construction of linear bases of free right-symmetric algebras is
given in \cite{Segal}.  Some further properties of this basis were
established in \cite{Dzhuma1}. Identities of right-symmetric
algebras were considered in \cite{Dzhuma2}. An analog of the
Magnus embedding theorem for right-symmetric algebras is proved in
\cite{KozyUmir}. The structure of universal multiplicative
enveloping algebras of free right-symmetric algebras is studied in
\cite{Kozy}.

In this paper we continue the study of free right-symmetric
algebras. We prove the Freiheitssatz for right-symmetric algebras
and the decidability of the word problem for right-symmetric
algebras with a single defining relation. These results generalize
the results of A.\,I.\,Shirshov mentioned above.

One of the fundamental results about free associative algebras is
the Bergman centralizer theorem ( see \cite{Berg}) which says that
the centralizer of any nonconstant element is a polynomial algebra
in a single variable. This theorem plays a crucial role in the
study of algorithmic and combinatorial questions. An analogue of
this result for free Poisson algebras over a field of
characteristic zero is proved in \cite{MakarU2}. In free Lie
algebras the centralizer is just the linear subspace generated by
an element. As we will see a similar result is true for the free
right-symmetric algebras.

It is well known (see, for example \cite{Cohn}) that if two
elements of a free associative algebra do not generate freely a
free subalgebra then they commute, and therefore by the Bergman
centralizer theorem are polynomials in a third variable. An
analogue of this result for free Poisson algebras is formulated in
\cite{MakarU2} and remains open. In the free right-symmetric
algebras we show that two nonconstant elements generate a free
right-symmetric algebra either of rank one or rank two.

The question about the freeness of subalgebras of free right
symmetric algebras was raised by A.\,S.\,Dzhumadildaev (oral
communication). Using the methods of \cite{Um11}, the first author
proved that the variety of right-symmetric algebras is not
Nielsen-Schreier. He also constructed a five generated non-free
subalgebra of a free right-symmetric algebra (see \cite{Kozy}).

We prove here that two generated subalgebras of free
right-symmetric algebras are free. This is rather surprising:
right-symmetric algebras are the first non Nielsen-Schreier
variety with this property known to us. For example two generated
subalgebras of polynomial algebras and free associative algebras
are not necessarily free.

It is well known \cite{Czer,Jung,Kulk,Makar} that the
automorphisms of polynomial algebras and free associative algebras
in two variables are tame. The automorphisms of free Poisson
algebras in two variables over a field of characteristic zero are
also tame \cite{MLTU}. We prove that the automorphisms of two
generated free right-symmetric algebras are tame.

It is also known \cite{Umi25,Um19} that
polynomial algebras (consequently, free Poisson algebras) and free associative algebras in three
variables in the case of characteristic zero have wild
automorphisms.
On the other hand, in 1964 P.Cohn \cite{Cohn2} proved that the automorphisms of a
free Lie algebra with a finite set of generators are tame.
The question about the existence of wild automorphisms of free right-symmetric algebras of larger ranks remains open.

This paper is organized as follows.
In Section 2 we prove some elementary properties of the linear basis constructed in \cite{Segal}.
In Section 3 we study algebras with a single defining relations. In Section 4 we study  subalgebras and automorphisms.

\bigskip

\section{Arithmetics of good words}

\hspace*{\parindent}

Let $X=\{x_1,x_2,\ldots,x_n\}$ be a finite alphabet. Denote by
$X^*$ the monoid of all nonassociative words on $X$. Denote by
$d(u)$ the degree function on $X^*$ such that $d(x_i)=1$ for all
$i$. Every nonassociative word $u$ of degree $\geq 2$ can be can
be uniquely written as $u=u_1u_2$, where $d(u_1), d(u_2)<d(u)$.

Put $x_1<x_2<\ldots <x_n$. Let $u$ and $v$ be arbitrary elements
of $X^*$. We say that $u<v$ if $d(u)<d(v)$. If $d(u)=d(v)\geq 2$,
$u=u_1u_2$, and $v=v_1v_2$, then $u<v$ if either $u_1<v_1$ or
$u_1=v_1$ and $u_2<v_2$.

A word is called {\em bad} if it contains a subword of the form
$(rs)t\in X^*$, where \\
 $d(r),d(s),d(t)\geq 1$ and $s>t$. A word is
called {\em good} if it is not bad. Denote by $W$ the set of all
good words in the alphabet $X$.

Let $A=A_n=RS\langle x_1,x_2,\ldots,x_n\rangle$ be the free
right-symmetric algebra in the variables $x_1,x_2,\ldots,x_n$ over
a field $k$. Every nonassociative word in the alphabet $X$
represents a certain element of $A$ and this representation gives
an embedding of $X^*$ into $A$. So for $u\in X^*$ we denote the
element of $A$ defined by $u$ by the same symbol.

According to \cite{Segal} the set of all good words form a linear basis of $A$:
every nonzero element $f$ of $A$ can be
uniquely represented as
\bes
f=\lambda_1 w_1+\lambda_2 w_2+\ldots +\lambda_m   w_m  ,
\ees
where $w_i\in W$, $0\neq \lambda_i\in k$ for all $i$ and $w_1>w_2>\ldots >w_m  $.

Denote by $\overline{f}$ the {\em leading word} $w_1$ of $f$.
The coefficient $\lambda_1$ is called the {\em leading coefficient} and
$\lambda_1w_1$ is called the leading term of $f$.
Usually we will assume that the leading coefficients
of the elements under consideration are equal to $1$.

\begin{lm}\label{l2}
Let $w\in X^*$ be an arbitrary nonassociative word.
Then $\overline{w}\leq w$ and the equality holds if and only if $w\in W$.
\end{lm}
\Proof  Assume that the statement of the lemma is true for all
nonassociative words $u$ such that $u<w$. If $w\in W$ then
$\overline{w}=w$.  Suppose that $w=uv$ is bad. If one of the words
$u$ or $v$ is bad then, by the assumption above, $w$ is a linear
combination of words $u'v'$ such that $u'v'<w$. Applying our
assumption to the words $u'v'$, we can say that $w$ is a linear
combination of good words $w'$ such that $w'<w$.

Suppose that both $u$ and $v$ are good. This means $u=u_1u_2$ and
$u_2>v$ since $w$ is bad. By (\ref{f1}), we have \bes
w=(u_1u_2)v=(u_1v)u_2+u_1(u_2v)-u_1(vu_2) \ees in the algebra $A$.
Note that $(u_1v)u_2,u_1(u_2v),u_1(vu_2)<w$. Then our assumption
gives the statement of the lemma again. $\Box$

\begin{rem}\label{r0}
$d(\overline{w}) = d(w)$ for any word $w \in X^*$.
\end{rem}
\noindent
Clear since all four words in (\ref{f1}) have the same
degree.

For every $f\in A$ denote by $R_f$ the operator of right
multiplication by $f$ acting on $A$, i.e., $uR_f=uf$ for all $u\in
A$. In particular, if $w,w_1,w_2,\ldots,w_m  \in X^*$ then
\bes
wR_{w_1}R_{w_2}\ldots R_{w_m  }= (\ldots ((ww_1)w_2)\ldots w_m  ).
\ees

\begin{lm}\label{l1}
A good word $w\in W$ can be uniquely represented as
\bee\label{f2}
w=x_iR_{w_1}R_{w_2}\ldots R_{w_m  },
\eee
where $w_j\in W$ for all $j$ and $w_1\leq w_2\leq \ldots \leq w_m  $.
\end{lm}
\Proof We prove it by induction on $d(w)$. If $w=uv$ then $u$ and
$v$ are good words and by induction hypothesis we can assume that
$u=x_iR_{u_1}R_{u_2}\ldots R_{u_s}$ where $u_j\in W$ and $u_1\leq
u_2\leq \ldots \leq u_s$. We have $u=u'u_s$ where
$u'=x_iR_{u_1}R_{u_2}\ldots R_{u_{s-1}}$. Note that $u_s\leq v$
since $w=uv=(u'u_s)v$ is good. Consequently,
$w=x_iR_{u_1}R_{u_2}\ldots R_{u_s}R_v$ and $u_1\leq u_2\leq \ldots
\leq u_s\leq v$. The base of induction when $d(w) = 1$ and
uniqueness are obvious. $\Box$

\begin{lm}\label{l3}
Let $u$ and $v$ be arbitrary good words and assume that
$u=x_iR_{u_1}R_{u_2}\ldots R_{u_m}$. Then
\bes
\overline{uv}=x_iR_{u_1}\ldots R_{u_s}R_vR_{u_{s+1}}\ldots
R_{u_m},
\ees
where $u_1\leq \ldots \leq u_s\leq v<u_{s+1}\leq
\ldots \leq u_m$.
\end{lm}
\Proof We assume that the statement of the lemma is true for all
pairs of good words $u'$ and $v'$ such that $u'v'<uv$. If $u_m
\leq v$ then $uv$ is also good and $\overline{uv}=uv$ satisfies
the statement of the lemma. Suppose that $u_m  >v$ and put
$w=x_iR_{u_1}R_{u_2}\ldots R_{u_{m-1}}$. Then, by (\ref{f1}), we
have \bee\label{f3} uv=(wu_m  )v=(wv)u_m  +w(u_m  v)-w(vu_m). \eee
By our assumption, $\overline{wv}=x_iR_{u_1}\ldots
R_{u_s}R_vR_{u_{s+1}}\ldots R_{u_{m-1}}$, where $u_1\leq \ldots
\leq u_s\leq v<u_{s+1}\leq \ldots \leq u_{m-1}$. Consequently,
\bes wv=\overline{wv}+\sum_i \alpha_i w_i, \ees where $w_i$ are
good and $\overline{wv}>w_i$ for all $i$. We have \bes (wv)u_m
=\overline{wv}u_m  +\sum_i \alpha_i w_iu_m  . \ees Note that $t
=\overline{wv}u_m$ is good and that $t>w_iu_m  $. Also $t>w(u_m
v),w(vu_m  )$ since $d(\overline{wv}) = d(w) + d(v)$ (see Remark
\ref{r0}). Lemma \ref{l2} gives  $\overline{(wv)u_m  } =t$. Then
(\ref{f3}) and Lemma \ref{l2} give that  $\overline{uv}=t$. $\Box$

\begin{co}\label{c0}
If $u, \ v, \ w \in W$ then $\overline{(uv)w}  = \overline{(uw)v}$.
\end{co}

\begin{lm}\label{l4}
Let $u$, $v$, and $w$ be arbitrary good words. If $u<v$ then
$\overline{wu}<\overline{wv}$ and $\overline{uw}<\overline{vw}$.
\end{lm}
\Proof We should consider only the case $d(u)=d(v)$, otherwise the
statement follows from the definition of order $\leq$. First we
prove that $\overline{wu}<\overline{wv}$ by induction on $d(w)$.
If $wv$ is good then $\overline{wv} = wv > wu \geq \overline{wu}$
(see Lemma \ref{l2}). If $wv$ is bad then $w=w_1w_2$ and $w_2 > v
> u$. In this case $\overline{wu}=\overline{w_1u}w_2$ and
$\overline{wv}=\overline{w_1v}w_2$ by Lemma \ref{l3}. By induction
on $d(w)$ we can assume that $\overline{w_1u}<\overline{w_1v}$.
Consequently, $\overline{wu}<\overline{wv}$.

Now we prove that $\overline{uw}<\overline{vw}$. If $vw$ is good
then Lemma \ref{l2} gives the claim. If $vw$ is bad then
$v=v_1v_2$ and $v_2>w$. By Lemma \ref{l3},
$\overline{vw}=\overline{v_1w}v_2$. Since $d(u) = d(v)$ we can
write $u = u_1 u_2$.

If $w > u_2$ then $\overline{uw}=(u_1u_2)w$. In this case $u_1u_2
< u_1 w \leq v_1 w$ and $u_1 u_2 = \overline{u_1u_2} <
\overline{u_1 w} \leq \overline{v_1 w}$ by the first claim of the
Lemma and by induction on $d(v)$. Therefore $(u_1u_2)w <
\overline{v_1w}v_2$ since $d((u_1u_2)w) = d(\overline{v_1w}v_2)$
and $u_1u_2 < \overline{v_1w}$.

If $w = u_2$ and $v_1 > u_1$ then inequalities above again give
$u_1u_2 < \overline{v_1w}$ and $(u_1u_2)w < \overline{v_1w}v_2$.

If $w = u_2$ and $v_1 = u_1$ then $v_2 > u_2$ and $\overline{uw} =
(v_1u_2)u_2 < (v_1u_2)v_2 = \overline{vw}$.

If  $u_2 > w$ then $\overline{uw}=\overline{u_1w}u_2$ by Lemma \ref{l3}.
If $u_1 < v_1$ then $\overline{u_1w}<\overline{v_1w}$ by induction and
$\overline{uw}=\overline{u_1w}u_2<\overline{v_1w}v_2=\overline{vw}$.
If $u_1 = v_1$ then $u_2<v_2$ and
$
\overline{uw}=\overline{u_1w}u_2=\overline{v_1w}u_2<\overline{v_1w}v_2=\overline{vw}.
$
The second claim of the lemma is proved. $\Box$

\begin{co}\label{c1}
If $f, g \in A$ then $\overline{fg}=\overline{\overline{f}\overline{g}}$.
\end{co}
\Proof If $w, \ x, \ y, \ z \in W$ and $w < x$, $y < z$ then $wy <
wz < xz$ and $\overline{wy} < \overline{wz} < \overline{xz}$ by
Lemma \ref{l4}. Since $f=\lambda_1 u_1+\lambda_2 u_2+\ldots
+\lambda_m u_m$ and $g = \mu_1 w_1+\mu_2 w_2+\ldots +\mu_{m'}
w_{m'} $ where $u_i, \ w_i\in W$, $\lambda_i, \ \mu_i \in k
\setminus 0$ for all $i$, $u_1 > u_2 > \ldots u_m$, and
$w_1>w_2>\ldots >w_{m'} $ we see that $\overline{u_1w_1} >
\overline{u_i w_j}$ if $(i, j) \neq (1, 1)$. $\Box$

\begin{lm}\label{l2'}
If $u$ and $v$ are non-empty good words and $u < v$ then
$\overline{uv} < \overline{vu}$.
\end{lm}
\Proof If $vu$ is a good word then $\overline{vu} = vu > uv \geq
\overline{uv}$. If $vu$ is a bad word then by Lemma \ref{l3}
$\overline{vu} = x_iR_{v_1} \ldots R_{v_s} R_u \ldots R_{v_m  }$
where $d(v_m  ) > 0$. So $\overline{vu} = wv_m  $ where $w =
x_iR_{v_1} \ldots R_{v_s} R_u \ldots R_{v_{k-1}}$ and $wv_m   >
uv$ since $d(uv) = d(\overline{vu})$ by Remark \ref{r0} and $d(w)
> d(u)$. $\Box$

\begin{co}\label{c1'}
If $f, g \in A \setminus k$ and $[f, g] = 0$ then $g = cf + \a$
where $c, \a \in k$. So $Z(f) = kf + k$ where $Z(f)$ is the
centralizer of $f$.
\end{co}
\Proof If $\overline{f} \neq \overline{g}$ then by the Lemma
$\overline{fg} \neq \overline{gf}$ and $[f, g] \neq 0$. So
$\overline{f} = \overline{g}$ and $\overline{g - cf} <
\overline{f}$ for some $c \in k$. Since $[f, g - cf] = 0$ we can
conclude that $g - cf = \a$. $\Box$

\newpage
\section{Algebras with a single defining relation}

\hspace*{\parindent}

In this section we prove that the word problem for right-symmetric
algebras with single relation is decidable and the Freiheitssatz.

Denote by $(f)$ the two-sided ideal of $A = RS\langle x_1, x_2,
\ldots, x_n\rangle$ generated by an $f \in A$. We would like to
find a linear basis of $(f)$.

Consider $B = RS\langle x_1, x_2, \ldots, x_n, y\rangle$ where $y$
is an additional variable.

Extend the order $<$ from $A$ to $B$ by assigning $d(y) = 1$ and
$y > x_n$ and define good words in the alphabet $\{x_1, \ldots,
x_n, y\}$ relative to this order exactly as it was done in Section
2. Denote this set by $W(B)$. For a non-associative word $w$ in
the alphabet $x_1, x_2, \ldots, x_n, y$ denote by $d_y(w)$ the
degree of $w$ relative to $y$.

Denote by $W_y$ the set of all good words $u \in W(B)$ with $d_y(u) = 1$.

\begin{lm}\label{l5}
If $\pi: B \mapsto A$ is a map defined by $\pi(x_i) = x_i$ and
$\pi(y) = f$ then the linear span of $\pi(W_y)$ is  $(f)$.
\end{lm}
\Proof Let $(y)$ be the ideal generated by $y$ in $B$. Then
$\pi((y)) = (f)$. Any element $Y \in (y)$ is a linear combination
of good words containing $y$. If we replace in these words  all
but one appearance of $y$ by $f$ then the image of the
corresponding linear combination of elements of $B$ will be
$\pi(Y)$. Since all modified words have degree 1 relative to $y$
they can be presented as linear combinations of elements of $W_y$.
$\Box$

Denote $\pi(W_y)$ by $E$.

\begin{lm}\label{l6}
There exists a subset $E'$ of $E$ which is a basis of $(f)$ such
that $\overline{g} \neq \overline{h}$ for $g \neq h \in E'$.
\end{lm}
\Proof We will lead induction on $g \in E$: we assume that the
subspace spanned by $E_g = \{h \in E | \overline{h} <
\overline{g}\}$ admits a basis consisting of elements of $E_g$
with different leading words. Since $E \subset X^*$ is a
well-ordered set, and if we take $g \in E$ with $\overline{g}$
minimal possible then the set $E_g$ is empty and satisfies the
Lemma.

In order to prove the Lemma we will check that if $g, \ h \in E$
and $\overline{g} = \overline{h}$ then $g - h$ is a linear
combination of elements $g_i$ from $E_g$.

Denote by $G, H \in W_y$ some elements for which $\pi(G) = g$ and
$\pi(H) = h$.

$G$ can be written either as $G = x_i R_{G_1} \ldots R_{G_l}$
where only $G_a$ contains $y$ or as $G = y R_{G_1} \ldots R_{G_l}$
where $d_y(G_i) = 0$  for $i = 1, \ldots, l$ and in both cases
$G_1 \leq G_2 \leq \ldots G_l$ (see Lemma \ref{l1}). Of course $H$
also can be written in one of these forms: $H = x_i R_{H_1} \ldots
R_{H_m}$ where only $H_b$ contains $y$ or $H = y R_{H_1} \ldots
R_{H_m}$ and $H_1 \leq H_2 \leq \ldots \leq H_m$.

If $G = y R_{G_1} \ldots R_{G_l}$ and $H = y R_{H_1} \ldots
R_{H_m}$ then $g = fR_{G_1} \ldots R_{G_l}$ and  \\
$h = f R_{H_1}
\ldots R_{H_m}$. Since $\overline{g} = \overline{h}$ Lemmas
\ref{l1} and \ref{l3} imply that $l = m$ and $G_i = H_i$. So $G =
H$ and $g - h = 0$.

If $G = y R_{G_1} \ldots R_{G_l}$ and $H =  x_i R_{H_1} \ldots
R_{H_m}$ then $g = fR_{G_1} \ldots R_{G_l}$ and\\  $h = x_i
R_{H_1} \ldots R_{\pi(H_b)} \ldots R_{H_m}$. Let $f = x_iR_{f_1}
\ldots R_{f_s}$. Since $\overline{\pi(H_b)} \geq f$ we should have
$\overline{\pi(H_b)} = G_j$ for some $j$. So $G_j = \pi(H_b) +
\delta$ where $\delta$ is a linear combination of good words which
are smaller than $G_j$. So
\bes
g = \pi(fR_{G_1} \ldots R_{G_{j-1}}
R_{H_b}R_{G_{j+1}} \ldots  R_{G_l} + yR_{G_1} \ldots R_{G_{j-1}}
R_{\delta}R_{G_{j+1}} \ldots  R_{G_l}).
\ees
 Since
\bes
\overline{fR_{G_1}
\ldots R_{G_{j-1}} R_{\delta}R_{G_{j+1}} \ldots  R_{G_l}} <
\overline{g}
\ees
 we can replace $G$ by $G' = fR_{G_1} \ldots
R_{G_{j-1}} R_{H_b}R_{G_{j+1}} \ldots  R_{G_l}$ and $g$ by
$\pi(G')$.

If $G = x_i R_{G_1} \ldots R_{G_l}$ and $H =  x_j R_{H_1} \ldots
R_{H_m}$ then $g =  x_i R_{G_1} \ldots R_{\pi(G_a)} \ldots
R_{G_l}$ and $h =  x_j R_{H_1} \ldots R_{\pi(H_b)} \ldots
R_{H_m}$. Since $\overline{g} = \overline{h}$ we should have $x_i
= x_j$ and $l = m$. Furthermore, either $\overline{\pi(G_a)} =
H_c$ and $\overline{\pi(H_b)} = G_d$ where $c \neq b$ and $d \neq
a$ or $\overline{\pi(G_a)} = \overline{\pi(H_b)}$.

If $\overline{\pi(G_a)} \neq \overline{\pi(H_b)}$ then
$G_d = \pi(H_b) + \delta$ where $\delta$ is a linear
combination of good words which are smaller than $G_d$.
Therefore
\bes
g = \pi(x_iR_{G_1} \ldots R_{G_{a-1}} R_{\pi(G_a)} R_{G_{a+1}} \ldots R_{G_{d-1}} R_{H_b}R_{G_{d +1}} \ldots  R_{G_l} +\\
x_iR_{G_1} \ldots R_{G_{a-1}} R_{G_a} R_{G_{a+1}} \ldots R_{G_{d-1}} R_{\delta}R_{G_{d+1}} \ldots  R_{G_l}).
\ees
Since
\bes
\overline{\pi(x_iR_{G_1} \ldots R_{G_{a-1}} R_{G_a}
R_{G_{a+1}} \ldots R_{G_{d-1}} R_{\delta}R_{G_{d+1}} \ldots
R_{G_l})} < \overline{g}
\ees
 we can replace $G$ by $G' = x_iR_{G_1}
\ldots R_{G_{a-1}} R_{\pi(G_a)} R_{G_{a+1}} \ldots R_{G_{d-1}}
R_{H_b}R_{G_{d +1}} \ldots  R_{G_l}$ and $g$ by $\pi(G')$. Here
$G_a' = H_b$.

Now consider the case when $\overline{\pi(G_a)} =
\overline{\pi(H_b)}$. Let $g_a = \pi(G_a)$, $h_b = \pi(H_b)$.
Since $\overline{g_a} < \overline{g}$ we can apply induction to
this pair and write $g_a = h_b + \delta$ where $\delta$ belongs to
the span of $E_{g_a}$. Therefore $\pi(G_a) = \pi(H_b + \Delta)$
where $\Delta$ belongs to the span of $W_y$ and $\pi(\Delta) \in $
span$(E_{g_a})$. So
\bes
G = x_iR_{G_1} \ldots R_{G_{a-1}} R_{H_b}
R_{G_{a+1}} \ldots  R_{G_l} + x_iR_{G_1} \ldots R_{G_{a-1}}
R_{\Delta} R_{G_{a+1}} \ldots  R_{G_l}.
\ees
 Since
\bes
\overline{\pi(x_iR_{G_1} \ldots R_{G_{a-1}} R_{\Delta}
R_{G_{a+1}} \ldots  R_{G_l})} < \overline{g}
\ees
 we can replace $G$
by $G' = x_iR_{G_1} \ldots R_{G_{a-1}} R_{H_b} R_{G_{a+1}} \ldots
R_{G_l}$.

It remains to consider the case when $G_a = H_b$. Since
$\overline{g} = \overline{h}$ the sets $G_1 \leq \ldots \leq
G_{a-1} \leq G_{a+1} \leq \ldots G_m$ and $H_1 \leq \ldots \leq
H_{b-1} \leq H_{b+1} \leq \ldots H_m$ should coincide. If $a = b$
then $G = H$ and the Lemma is proved. Assume that $a > b$. Then
$G_a > G_{a-1} = H_a > H_b$ ($G_a \neq G_{a-1}$ since
$d_y(G_{a-1}) = 0$).
But $G_m = H_b$ which proves the Lemma. $\Box$

The undecidability of the  word problem for Lie algebras was
proved by L.\,A.\,Bokut' \cite{Bokut}. An explicit example of a
finitely presented Lie algebra with undecidable word problem was
constructed by G.\,P.\,Kukin \cite{Kukin}. If $L$ is a finitely
presented Lie algebra then the universal right-symmetric
enveloping algebra $A(L)$  of $L$ is also finitely presented
algebra with the same set of generators and defining relations
\cite{Segal}.  Consequently, the word problem for right-symmetric
algebras is also undecidable. On the other hand, A.\,I.\,Shirshov
\cite{Shir2} proved the decidability of the word problem for Lie
algebras with a single defining relation. In the case of
right-symmetric algebras we have the next result.

\begin{theor}\label{t1}
The word problem for right-symmetric algebras with a single
defining relation is decidable.
\end{theor}

\Proof
Let us fix a basis
$E' = \{e_1, e_2, \ldots, e_m, \ldots \}$ with $\overline{e_1} < \overline{e_2} \ldots < \overline{e_m} < \ldots$
existence of which is guaranteed by Lemma \ref{l6}.
Take an $h \in A$. If $h \in (f)$ then $h = \sum_{i=1}^a \lambda_i e_i$.

Let $d_i = d(\overline{e_i})$. Since $\overline{e_1} <
\overline{e_2} \ldots < \overline{e_m} < \ldots$, the sequence of
natural numbers $d_1, d_2, \ldots, d_m, \ldots$ is a nondecreasing
sequence and for a natural number $l$ there is only finitely many
elements of this sequence equal to $l$. So if $h \in (f)$ then $h$
belongs to a finite-dimensional space $V_d = $ span$\{ e_i \in E'|
d_i \leq d(\overline{h})\}$. Therefore we can effectively
determine whether $h$ is in $(f)$. $\Box$

\begin{rem}\label{r1}
Since $d_i \geq d(\overline{f})$ an ideal $(f)$ does not contain
elements $h$ with $d(\overline{h}) < d(\overline{f})$.
\end{rem}

As we mentioned in the Introduction, the Freiheitssatz for Lie
algebras was proved by A.\,I.\,Shirshov \cite{Shir2} and
L.\,Makar-Limanov \cite{Makar2} proved the Freiheitssatz for free
associative algebras in the case of characteristic zero. The
analogue of these results for right-symmetric algebras is also
true.

\begin{theor}\label{t2}{\bf (Freiheitssatz)}
If $f \in RS\langle x_1,x_2,\ldots,x_n\rangle$ and $f\notin
RS\langle x_1,x_2,\ldots,x_{n-1}\rangle$, then $(f)\cap RS\langle
x_1,x_2,\ldots,x_{n-1}\rangle=0$.
\end{theor}
\Proof Let $h \in (f) \cap RS\langle
x_1,x_2,\ldots,x_{n-1}\rangle$. For a $w \in W$ consider an
endomorphism $\rho_w$ of $A$ to $A$ given by $\rho_w(x_i) = x_i$
if $i < n$ and $\rho_w(x_n) = x_n w$. It is clear that $h \in
(\rho_w(f)) \cap RS\langle x_1,x_2,\ldots,x_{n-1}\rangle$ for any
$w$. If $h \neq 0$ take a $w$ with $d(w) > d(\overline{h})$. Then
$d(\overline{\rho_w(f)}) > d(\overline{h})$ since $f$ contains
$x_n$. Therefore by Remark \ref{r1}, $(\rho_w(f)) \not\ni h$.
$\Box$

\section{Two generated subalgebras and automorphisms}

\hspace*{\parindent}

As in the preceding sections $A=RS\langle
x_1,x_2,\ldots,x_n\rangle$ and $W$ is the set of all good words in
the alphabet $X=\{x_1,x_2,\ldots,x_n\}$.

Recall that a set of elements $s_1,s_2,\ldots,s_m$ of a polynomial
algebra is called {\em algebraically independent} if the
subalgebra generated by $s_1,s_2,\ldots,s_m $ is a polynomial
algebra in the variables $s_1,s_2,\ldots,s_m $. We will use
analogous terminology in the case of other free algebras. A set of
elements $s_1,s_2,\ldots,s_m$ of a free right-symmetric (Lie,
associative, or Poisson) algebra is called {\em free} if the
subalgebra generated by $s_1,s_2,\ldots,s_m$ is free  and
$s_1,s_2,\ldots,s_m$ is a free set of generators of this
subalgebra.

\begin{lm}\label{l10}
Let $z \in A$ be an arbitrary element which does not belong to the
filed $k$. Then the subalgebra of $A$ which is generated by $z$ is
isomorphic to $RS\langle x\rangle$
\end{lm}
\Proof Assume that it is not the case. Then there exists a
non-zero element $p(x) \in RS\langle x\rangle$ such that $p(z) =
0$. Since $p(x) = \sum \lambda_i w_i$ where $\lambda_i \in k$ and
$w_i$ are good words in alphabet $\{x\}$ we should have two
different words $w_a(x)$ and $w_b(x)$ for which $\overline{w_a(z)}
= \overline{w_b(z)}$. Let us assume that the pair $w_a(x), \
w_b(x)$ is a minimal pair with this property. We can write $w_a(x)
= x R_{u_1} \ldots R_{u_l}$ and $w_b(x) = x R_{v_1} \ldots
R_{v_m}$ according to Lemma \ref{l1}. So $\overline{z R_{u_1}
\ldots R_{u_l}} = \overline{z R_{v_1} \ldots R_{v_m}}$ and Lemma
\ref{l3} implies that $l = m$ and for each $u_i$ there is a $v_j$
such that $\overline{u_i(z)} = \overline{v_j(z)}$. Since $d(u_i) <
d(w_a)$ and $d(v_j) < d(w_b)$ we can conclude that $u_i(x) =
w_j(x)$. But then good words $w_a$ and $w_b$ are equal. $\Box$

\begin{theor}\label{t3}
Two generated subalgebras of free right-symmetric algebras are free.
\end{theor}
\Proof Let $f_1, f_2 \in A$. Assume that $f_1$ and $f_2$ are
dependent, i. e. there exists a non-zero element $p(y_1, y_2) \in
RS\langle y_1,y_2\rangle$ such that $p(f_1, f_2) = 0$. If
$\overline{f_1} = \overline{f_2}$ we can chose $\lambda \in k$ so
that $\overline{f_2 - \lambda f_1} < \overline{f_1}$ and replace
the pair $f_1, \ f_2$ with a dependent pair $f_1, \ f_2 - \lambda
f_1$. Hence we assume that $\overline{f_1} > \overline{f_2}$.

Let $\pi$ be a homomorphism of $RS\langle y_1,y_2\rangle$ into
$RS\langle x_1,x_2,\ldots,x_n\rangle$ given by $\pi(y_1) = f_1$,
$\pi(y_2) = f_2$.

We can write $p = \sum \lambda_i w_i$ where $w_i \in W(y_1, y_2)$
and $\lambda_i \in k$. Since $\pi(p) = 0$ we see that there should
be a pair of different words $w_a$ and $w_b$ in this sum for which
$\overline{\pi(w_a)} = \overline{\pi(w_b)}$.

Let us assume that the pair $w_a, \ w_b$ is a minimal pair with
this property. We can write $w_a = y_1 R_{u_1} \ldots R_{u_l}$
renaming $y_1$ and $y_2$ if necessary. Let $w_b = y_1 R_{v_1}
\ldots R_{v_m}$. Since $\overline{\pi(w_1)} = \overline{\pi(w_2)}$
we have $\overline{f_1 R_{\pi(u_1)} \ldots R_{\pi(u_l)}} =
\overline{f_1R_{\pi(v_1)} \ldots R_{\pi(v_m)}}$. Then by Lemmas
\ref{l1} and \ref{l3} for each $u_i$ there is a $v_j$ such that
$\overline{\pi(u_i)} = \overline{\pi(v_j)}$. Since $w_a, \ w_b$ is
a minimal pair and $d(u_i) < d(w_a)$, $d(v_j) < d(w_b)$ we can
conclude that $u_i = v_j$. Since $u_1 \leq \ldots \leq u_l$ and
$v_1 \leq \ldots \leq v_m$ we should have $w_a = w_b$ contrary to
our assumption. So $w_b(y_1, y_2) = y_2 R_{v_1} \ldots R_{v_m}$.
Therefore $\overline{f_1 R_{\pi(u_1)} \ldots R_{\pi(u_l)}} =
\overline{f_2R_{\pi(v_1)} \ldots R_{\pi(v_m)}}$.

If $\overline{\pi(u_i)} = \overline{\pi(v_j)}$ for some pair $i, \
j$ then by Lemma \ref{l3} $\pi(w_a') = \pi(w_b')$ where $w_a'$ is
$w_a$ with omitted $R_{u_i}$ and $w_b'$ is $w_b$ with omitted
$R_{v_j}$. Since $d(w_a') < d(w_a)$ and $d(w_b') < d(w_b)$ it
would imply that $w_a' = w_b'$. But this impossible since $w_a' =
y_1 R_{u_1} \ldots R_{u_l}$ and $w_b' = y_2 R_{v_1} \ldots
R_{v_m}$. Therefore $\overline{\pi(u_i)} \neq \overline{\pi(v_j)}$
for any pair $i, \ j$. Since $\overline{f_1 R_{\pi(u_1)} \ldots
R_{\pi(u_l)}} = \overline{f_2R_{\pi(v_1)} \ldots R_{\pi(v_m)}}$
Lemma s \ref{l1} and \ref{l3} imply that $\overline{f_1} =
\overline{x_i R_{\pi(v_1)} \ldots R_{\pi(v_m)} R_{g_1} \ldots
R_{g_m}}$ and $\overline{f_2} = \overline{x_i R_{\pi(u_1)} \ldots
R_{\pi(u_l)} R_{g_1} \ldots R_{g_m}}$ where $g_i \in A$.

Since $d(\overline{f_1}) \geq d(\overline{f_2})$ we see that
$d(\pi(u_i)) \geq d(f_2)$ if $d(u_i) > 0$. Therefore $f_2 =
R_{g_1} \ldots R_{g_m}$. If $v_i$ contains $y_1$ then $d(\pi(v_i))
\geq d(f_1)$ which is also impossible. So $v_i \in RS\langle
y_2\rangle$ and $\overline{f_1} \in RS\langle
\overline{f_2}\rangle$. But then with the right choice of $q$ we
can replace the pair $f_1, \ f_2$ by the pair $f_1 - q(f_2), \
f_2$ where $\overline{f_1 - q(f_2)} < \overline{f_1}$. We can
conclude by induction that the subalgebra of $A$ generated by
$f_1, \ f_2$ is $RS\langle g\rangle$ for some element $g \in A$.
$\Box$

A pair of elements $f$ and $g$ of the algebra $A$ is called {\em
reducible} if there exists a good word $s$ in the variable $x$
such that $s(\overline{f})=\overline{g}$ or
$s(\overline{f})=\overline{g}$.  A pair $f,g$ is called {\em
reduced} if it is not reducible.

Consider a subalgebra $S$ of $A$ generated by two nonzero elements
$f$ and $g$. If the pair $f$ and $g$ is reduced and both $f$ and
$g$ are not in $k$ then by Theorem \ref{t3} they generate a free
right-symmetric subalgebra of $A$ of rank two and $f$ and $g$ are
free generators of this subalgebra.

Recall that an automorphism $\phi$ of a free right-symmetric
algebra $A$ generated by $\{ x_1,x_2,\ldots,x_n\}$ is called {\em
elementary} if $\phi(x_j) = x_j$ for any $j \neq i$ and $\phi(x_i)
= \a x_i + f$ where $f\in RS\langle
x_1,\ldots,x_{i-1},x_{i+1},\ldots,x_n\rangle$. Automorphisms which
can be expressed as a composition of elementary automorphisms are
called {\em tame}. Non-tame automorphisms are called {\em wild}.

Denote by $\phi = (f_1,f_2,\ldots,f_m )$ an automorphism $\phi$ of
$A$ such that $\phi(x_i)=f_i,\, 1\leq i\leq n$. It is well known
(see, for example \cite{Cohn}) that $\phi$ is tame if and only if
there exists a sequence of elementary transformations such that
\bes
(f_1,f_2,\ldots,f_m)=\psi_r\rightarrow\psi_{r-1}\rightarrow\psi_{r-2}\rightarrow\ldots
\rightarrow\psi_0= (x_1,x_2,\ldots,x_n). \ees

\begin{theor}\label{t4}
Automorphisms of two generated free right-symmetric algebras are
tame.
\end{theor}
\Proof Let $\varphi=(f_1, f_2)$ be an automorphism of $A_2 =
RS\langle x_1,x_2\rangle$. If the pair $f_1, \ f_2$ is reducible
then, using an appropriate elementary reduction $\phi_{1,q}(f_1) =
f_1$, $\phi_{1,q}(f_2) = f_2 - q(f_1)$ or $\phi_{2,q}(f_1) = f_1 -
q(f_2)$, $\phi_{2,q}(f_2) = f_2$ we can decrees
$d(\overline{f_1f_2})$. So after a finite number of reductions we
obtain a pair $f_1', \ f_2'$ which is reduced and still generate
$A_2$. So $x_1 = p_1(f_1', f_2')$ and $x_2 = p_2(f_1', f_2')$.
Since different good words in $f_1', \ f_2'$ have different
leading words in $A_2$ (see the proof of Theorem \ref{t3}) we see
that $x_1 = \overline{x_1} = \overline{w_1(f_1', f_2')}$. So $x_1
= \overline{f_i' R_{u_1} \ldots R_{u_l}}$ and $x_2 =
\overline{f_j' R_{v_1} \ldots R_{v_m}}$. This is possible only if
$x_1 = \overline{f_i'}$ and $x_2 = \overline{f_j'}$. So $f_i' =
\lambda_1 x_1 + \mu_1$ and $f_j' = \lambda_2 x_2 + \mu_2 x_1 +
\nu_2$ where $\lambda_1, \ \lambda_2, \ \mu_1, \ \mu_2, \ \nu_2
\in k$, $\lambda_1, \ \lambda_2 \neq 0$ and $\varphi$ is tame.
$\Box$

Now we want to formulate some open questions closely related to obtained results.
\begin{prob}\label{pr1}
Are three generated subalgebras of free right-symmetric algebras free?
\end{prob}
Note that a five generated non-free subalgebra of a free
right-symmetric algebra was given in \cite{Kozy}.

It is well known that the algebraic dependence of a finite set of
elements of a polynomial algebra is algorithmically recognizable.
There exists an algorithm which decides whether a finite set of
elements in a free Lie algebra is free (see, for example
\cite{Mikhalev}). It is also known that the freeness of a finite
set of elements is algorithmically unrecognizable for free
associative algebras \cite{UUU}.
\begin{prob}\label{pr2}
Is the freeness of a finite set of elements of a free
right-symmetric algebra algorithmically recognizable?
\end{prob}
\begin{prob}\label{pr3}
Are the automorphisms of finitely generated free right-symmetric algebras tame?
\end{prob}

\bigskip

\begin{center}
{\bf\large Acknowledgments}
\end{center}

\hspace*{\parindent}

The third author wishes to thank the Department of Mathematics of
Wayne State University in Detroit for the support while he was
working on this project.

\end{document}